\catcode`\@=11

\magnification 1200

\pretolerance=500 \tolerance=1000 \brokenpenalty=5000

\hsize=134mm \vsize=204mm
\hoffset=0mm \voffset=0mm

\parindent 6mm

\parskip 1,2mm

\baselineskip=12,8pt

\hfuzz 2mm

\newcount\equacount
\def\ifundefined#1{\expandafter\ifx\csname#1\endcsname\relax}
\def\equadef#1{\global\advance\equacount by 1%
  \global\expandafter\edef\csname#1\endcsname{\the\equacount}%
  \the\equacount}
\def\equaref#1{\expandafter\csname#1\endcsname}

\newif\ifpagetitre        \pagetitretrue
\newtoks\hautpagetitre    \hautpagetitre={\hfil}
\newtoks\baspagetitre     \baspagetitre={\hfil}
\newtoks\auteurcourant    \auteurcourant={\hfil}
\newtoks\titrecourant     \titrecourant={\hfil}
\newtoks\hautpagegauche   \newtoks\hautpagedroite
\hautpagegauche={\hfil\tensl\the\auteurcourant\hfil}
\hautpagedroite={\hfil\tensl\the\titrecourant\hfil}

\newtoks\baspagegauche
\baspagegauche={\hfil\tenrm\folio\hfil}
\newtoks\baspagedroite
\baspagedroite={\hfil\tenrm\folio\hfil}

\headline={\ifpagetitre\the\hautpagetitre
\else\ifodd\pageno\the\hautpagedroite
\else\the\hautpagegauche\fi\fi}

\footline={\ifpagetitre\the\baspagetitre
\global\pagetitrefalse
\else\ifodd\pageno\the\baspagedroite
\else\the\baspagegauche\fi\fi}

\font\bbb =msbm10
\def \rr {{\hbox {\bbb R}}}

\font\bbbb =msbm7
\def \rrr {{\hbox {\bbbb R}}}

\font\bbb =msbm10

\font\bbbb =msbm7

\font\bbb =msbm10

\font\bbbb =msbm7
\def \nnn {{\hbox {\bbbb N}}}

\font\bbb =msbm10

\font\bbbb =msbm7

\font\bbb =msbm10

\font\bbb =msbm10

\font\bbb =msbm10
\def \ss {{\hbox {\bbb S}}}

\font\bbbb =msbm7
\def \sss {{\hbox {\bbbb S}}}

\def\ref#1|#2|
       {{\leftskip=6mm\noindent
           \hangindent=5mm\hangafter=1
           \llap{[#1]}\hskip 5mm{#2}.\par}}

\def\pe#1#2 {\pc#1#2|}
\def\pc#1#2|{{\tenrm#1\sevenrm#2}}
\def\pd#1 {{\pc#1|}\ }
\def\bfpc#1#2|{{\tenbf#1\sevenbf#2}}
\def\bfpd#1 {{\bfpc#1|}\ }

\def\cqfd{\hfill\penalty 500\raise 0pt\hbox{\vrule\vbox to 6.5pt{
\hrule width 5.8pt\vfill\hrule}\vrule}\par}

\def \P {{(P_t)}_{t\geq 0}} 
\def \eu {{\rm e}}
\def \L {{\rm L}}

\centerline { }

\vskip 10mm

\centerline { \bf HYPERCONTRACTIVE MEASURES,}
\vskip 1mm
\centerline { \bf TALAGRAND'S INEQUALITY, AND INFLUENCES}

\vskip 0,5cm

\font\pc=cmcsc10 \rm 

\centerline {{\pc D. Cordero-Erausquin, M. Ledoux}}

\vskip 0,1cm

\centerline {\it University of Paris 6 and University of Toulouse, France}

\vskip 1cm

\parindent 1,2cm

{\narrower 
 \parindent 6mm

Abstract. -- {\it We survey several Talagrand type inequalities and their application
to influences with the tool of hypercontractivity for both discrete and continuous,
and product and non-product models. The approach covers
similarly by a simple interpolation the framework of geometric influences recently
developed by N. Keller, E. Mossel and ${\hbox {A. Sen}}$.
Geometric Brascamp-Lieb decompositions are also considered in this context.}
\par }

\parindent 6mm

\parskip 1,8mm

\vskip 8mm \goodbreak

{\bf 1. Introduction}

\vskip 3mm

In the famous paper [T], M. Talagrand showed that for every function $f$ on the discrete cube
$ X = \{-1, +1\}^N$ equipped with the uniform probability measure $\mu $,
$$ {\rm Var}_\mu (f) = \int_X  f^2 d\mu  - \bigg ( \int_X f d\mu \bigg)^2
      \leq  C \sum_{i=1}^N {  {\| D_i f\| }_2^2  
          \over 1+ \log \big ( {\| D_i f\| }_2 / {\| D_i f\| }_1 \big ) }  \eqno (1) $$
for some numerical constant $C \geq 1$,          
where ${\| \cdot \| }_p$ denote the norms in $\L^p(\mu )$, $1 \leq p \leq \infty$, and
for every $i = 1, \ldots, n$ and every $ x = (x_1, \ldots, x_N) \in  \{-1, +1\}^N$,
$$ D_i f(x) = f( \tau _i x) - f (x)  \eqno (2) $$
with $\tau _i x = (x_1, \ldots, x_{i-1}, -x_i, x_{i+1}, \ldots , x_N)$.
Up to the numerical
constant, this inequality improves upon the classical spectral gap inequality (see below)
$$ {\rm Var}_\mu (f) \leq  {1 \over 4} \sum_{i=1}^N   {\| D_i f\| }_2^2 \, . \eqno (3) $$

The proof of (1) is based on an hypercontractivity estimate known as the Bonami-Beckner
inequality [Bo], [Be] (see below).
Inequality (1) was actually deviced to recover (and extend) a famous result of J. Kahn, G. Kalai
and N. Linial [K-K-L] about influences on the cube. Namely, applying (1) to the Boolean
function $f = {\bf 1}_A $ for some set $A \subset \{-1, +1\}^N$, it follows that
$$ \mu (A) \big ( 1 - \mu  (A) \big ) 
   \leq  C \sum_{i=1}^N{ 2I_i(A) \over  1 + \log \big (1/ \sqrt {2I_i(A)} \, \big )} 
     \eqno (4) $$
where, for each $i = 1, \ldots, N$,
$$ I_i (A) = \mu \big ( \{ x \in A, \tau _i x \notin A \}  \big )$$
is the so-called influence of the $i$-th coordinate on the set $A$
(noticing that $\| D_i {\bf 1}_A \|^p _p = 2 I_i(A)$ for every $p \geq 1$).
In particular, for a set $A$ with $ \mu (A) = a$, there is a coordinate $i$, $1 \leq i\leq N$, such that
$$ I_i (A) \geq    {a(1-a) \over 8CN} \, \log \Big ( {N \over a(1-a)} \Big)
         \geq  {a(1-a) \log N  \over 8C N}  \eqno (5) $$
which is the main result of [K-K-L]. (To deduce (5) from (4), assume for example that
$I_i(A) \leq \big ( {a(1-a) \over N} \big) ^{1/2}$ for every $i=1, \ldots , N$,
since if not the result holds. Then, from (4), there exists $i$, $1 \leq  i \leq N$, such that
$$ {a(1-a) \over CN} \leq { 2I_i(A) \over  1 + \log \big (1/\sqrt { 2I_i(A)} \, \big )} 
     \leq { 8 I_i(A) \over  4 + \log ( N / 4 a(1-a) ) } $$
which yields (5)). Note that (5) remarkably improves by a (optimal) factor $\log N$ what would
follow from the spectral gap inequality (3) applied to $ f = {\bf 1}_A$.
The numerical constants like $C$ throughout this text are not sharp.

The aim of this note is to amplify the hypercontractive proof of Talagrand's original inequality
(1) to various settings, including non-product spaces and continuous variables, and
in particular to address versions suitable to geometric influences.
It is part of the folklore indeed (cf. e.g. [B-H]) that an inequality similar to (1), with the same
hypercontractive proof, holds for the standard Gaussian measure $\mu $
on $\rr^N$ (viewed as a product measure of one-dimensional
factors), that is, for every smooth enough function $f$ on $\rr^N$ and some constant $C>0$,
$$ {\rm Var}_\mu (f)  \leq  C \sum_{i=1}^N {  {\| \partial_i f\| }_2^2  
          \over 1+ \log ( {\| \partial_i f\| }_2 / {\| \partial_i f\| }_1) } \, .   \eqno (6) $$
(A proof will be given in Section 2 below.)          
However, the significance of the latter for influences is not clear, since its application
to characteristic functions is not immediate (and requires notions of capacities). Recently,
N. Keller, E. Mossel and A. Sen [K-M-S] introduced a notion of geometric influence of a Borel set
$A$ in $\rr^N$ with respect to a measure $\mu $ (such as the Gaussian measure) simply
as $ {\| \partial_i f \| }_1$  for some smooth approximation $f$ of $ {\bf 1}_A$,
and proved for it the analogue of (5) (with $\sqrt {\log N}$ instead of $\log N$) for the
standard Gaussian measure on $\rr^N$. It is therefore
of interest to seek for suitable versions of Talagrand's inequality involving only $\L^1$-norms
$ {\| \partial_i f \| }_1$ of the partial derivatives. While the authors of [K-M-S] use isoperimetric
properties, we show here how the common hypercontractive tool 
together with a simple interpolation argument may be 
developed similarly to reach the same conclusion. In particular, for the
standard  Gaussian measure $\mu $ on $\rr^N$, we will see that for every smooth
enough function $f $ on $\rr^N$ such that $|f| \leq 1$,
$$ {\rm Var}_\mu (f) \leq  C \sum_{i=1}^N   
       { {\| \partial _i f\|} _1 \big (1 + {\| \partial   _i f\|} _1  \big)
           \over  \big [ 1 +    \log^+ \big (1/  {\| \partial _i f\|}_1 \big ) \big ]^{1/2}  }  \, .  \eqno (7) $$
Applied to $f = {\bf 1}_A$, this inequality indeed ensures the existence of a coordinate $i$,
$1 \leq i\leq N$, such that the geometric influence of $A$ along $i$ is
at least of the order of $ {\sqrt {\log N} \over  N} $,
that is one of the main conclusions of [K-M-S] (where it is shown moreover that
the bound is sharp). In this continuous setting,
the hypercontractive approach yields more general examples of measures with such an
influence property in the range between exponential and Gaussian
for which only a logarithmic Sobolev type inequality is needed
while [K-M-S] required an isoperimetric inequality for the individual measures $\mu_i$.

This note is divided into two main parts. In the first one, we present Talagrand type inequalities
for various models, from the discrete cube to Gaussian and more general product measures,
by the general principle of hypercontractivity of Markov semigroups.
The method of proof, originating in Talagrand's work, has been used recently by R.~O'Donnell and K. Wimmer
[OD-W1], [OD-W2] to investigate non-product models such as random walks
on some graphs which enter the general presentation below. Actually, most
of the Talagrand inequalities we present in the discrete setting are already contained 
in the work by R. O'Donnell and K. Wimmer.
It is worth mentioning that an approach
to the Talagrand inequality (1) rather based on the logarithmic Sobolev inequality
was deviced in [Ros] and [F-S] a few years ago.
The abstract semigroup approach applies in the same way on the sphere
along the decomposition of the Laplacian.
Geometric Brascamp-Lieb decompositions within this setting are also discussed.
In the second part, we address our new version (7) of Talagrand's inequality towards
geometric influences and the recent results of [K-M-S] by a further
interpolation step on the hypercontractive proof.

In the last part of this introduction, we describe a convenient framework in order to develop
hypercontractive proofs of Talagrand type inequalities. While of some abstract flavor, the setting
easily covers two main concrete instances, probability measures on finite state spaces
(as invariant measures of some Markov kernels) and continuous probability measures
of the form $d\mu (x) = \eu^{-V(x)} dx$ on the Borel sets of $\rr^n$ where $V$ is some (smooth)
potential (as invariant measures of the associated diffusion operators
$ {\Delta - \nabla V \cdot \nabla }$). We refer for the material below to the general
references [Ba], [D-SC], [Roy], [Aal], [B-G-L]...

Let $\mu $ be a probability measure on a measurable space $(X, {\cal A})$. For a function
$f : X \to \rr$ in $\L^2(\mu )$, define its variance with respect to $\mu $ by
$$ {\rm Var}_\mu  (f) = \int_X f^2 d\mu  - \bigg ( \int_X f d\mu  \bigg)^2 .$$
Similarly, whenever $f >0$, define its entropy by
$$ {\rm Ent}_\mu  (f) = \int_X f \log f d \mu  - \int_X f d\mu \log \bigg ( \int_X f d\mu \bigg)$$
provided it is well-defined. The $\L^p(\mu )$-norms, $1 \leq p \leq \infty$, will be denoted by
${\| \cdot \| }_p$.

Let then $\P$ be a Markov semigroup with generator $\L $ acting on a suitable
class of functions on $(X, {\cal A})$. Assume that $\P $ and $\L$ have an invariant,
reversible and ergodic probability measure $\mu $. This ensures that the operators $P_t$ are
contractions in all $\L^p(\mu )$-spaces, $1 \leq  p \leq \infty$. The Dirichlet form
associated to the couple $(\L , \mu)$ is then defined, on functions $f, g$ of the Dirichlet
domain, as
$$ {\cal E} (f, g) = \int_X f (-\L g) d\mu . $$

Within this framework, the first example of interest is the case
of a Markov kernel $K$ on a finite state space $X$ with invariant
($ \sum_{x \in X} K(x,y) \mu (x) = \mu (y)$, $ x \in X$) and reversible
($K(x,y) \mu (x) = K(y,x) \mu (y)$,  $x,y \in X$) probability measure $\mu $.
The Markov operator $ \L = K- {\rm Id} $ generates the semigroup 
of operators $P_t = \eu^{t \L}$, $t\geq 0$, and defines the Dirichlet form
$$ {\cal E} (f,g)   = \int _X f (-\L g)  d\mu 
       = {1\over 2} \sum_{x,y \in X} \big [ f(x) - f(y) \big ] \big [g(x) - g(y) \big ] K(x,y) \mu (x) $$
on functions $f, g : X \to \rr$.       
The second class of examples is the case of $X= \rr^n$ equipped with its Borel $\sigma $-field.
Letting $V : \rr^n \to \rr$ be such that $\int_{\rrr^n} \eu^{-V(x)} dx = 1$,
under mild smoothness and growth conditions on the potential $V$, the
second order operator $\L = \Delta - \nabla V \cdot \nabla $ admits
$d\mu (x) = \eu^{-V(x)} dx$ as symmetric and invariant probability measure. The operator $\L $
generates the Markov semigroup of operators $\P$ and
defines by integration by parts the Dirichlet form
$$ {\cal E} (f,g) =  \int_{\rrr^n} f (- \L g)  d \mu  = \int_{\rrr^n} \nabla f \cdot \nabla g  \, d\mu $$
for smooth functions $f,g$ on $\rr^n$.

Given such a couple $(\L, \mu )$, it is said to satisfy a spectral gap, of
Poincar\'e, inequality if there is a constant $\lambda >0$ such that for all functions
$f$ of the Dirichlet domain,
$$ \lambda \, {\rm Var}_\mu (f) \leq {\cal E} (f,f) .  \eqno (8) $$
Similarly, it satisfies a logarithmic Sobolev inequality if there is a constant $\rho  >0$ such
that for all functions $f$ of the Dirichlet domain,
$$ \rho \, {\rm Ent}_\mu (f^2) \leq  2 \, {\cal E}(f,f). \eqno (9) $$
One speaks of the spectral gap constant (of $(\L, \mu )$) as the best $\lambda >0$ for which (8) holds,
and of the logarithmic Sobolev constant (of $(\L, \mu )$) as the best $\rho >0$ for which (9) holds.
We still use $\lambda $ and $\rho $ for these constants.
It is classical that $\rho \leq \lambda $.

Both the spectral gap and logarithmic Sobolev inequalities translate equivalently on the
associated semigroup $\P$. Namely, the spectral gap inequality (8) is
equivalent to saying that
$$ {\| P_t f \| }_2 \leq \eu^{-\lambda t}  \, {\| f \| }_2   $$
for every $t \geq 0$ and every mean zero function $f$ in $\L^2(\mu )$. Equivalently
for the further purposes, for every $f \in \L^2(\mu )$ and every $t > 0$,
$$ {\rm Var}_\mu (f) \leq  {1 \over 1 - \eu^{-\lambda t}} \,
         \big [ {\| f\| }_2^2 - {\| P_t f\| }_2^2 \big ] . \eqno (10) $$
On the other hand,
the logarithmic Sobolev inequality gives rise to hypercontractivity which is a smoothing property
of the semigroup. Precisely, the logarithmic Sobolev inequality (9) is equivalent to saying that,
whenever $p \geq 1 + \eu^{-2\rho t}$, for all functions $f$ in $\L^p(\mu )$,
$$ {\| P_t f \| }_2\leq  {\|  f \| }_p  .  \eqno (11) $$
For simplicity, we say below that a probability measure $\mu $ in this context is
hypercontractive with constant $\rho $. 

A standard operation on Markov operators is the product operation. Let
$(\L_1, \mu _1)$ and $(\L_2, \mu _2)$ be Markov operators on respective spaces
$X_1$ and $X_2$. Then
$$ \L = \L_1 \otimes {\rm Id} + {\rm Id} \otimes \L_2 $$
is a Markov operator on the product space $X_1 \times X_2$ equipped with the
product probability measure $\mu _1 \otimes \mu _2$. The product semigroup
$\P$ is similarly obtained as the tensor product $P_t = P_t^1 \otimes P_t^2$
of the semigroups on each factor. For the product Dirichlet form, the spectral
gap and logarithmic Sobolev constants are stable in the sense that, with the obvious notation,
$\lambda = \min (\lambda _1, \lambda _2)$ and $\rho = \min (\rho _1, \rho _2)$.
This basic stability by products will allow for constants independent of the dimension
in the Talagrand type inequalities under investigation. 
For the clarity of the exposition, we will not mix below products of continuous and
discrete spaces, although this may easily be considered.

Let us illustrate the preceding definitions and properties on two basic examples. Consider first
the two-point space $ X = \{-1, +1\}$ with the measure $\mu  = p \delta_{+1} + q \delta_{-1}$,
$p \in [0,1]$, $p+q=1$, and the Markov kernel $K(x,y ) = \mu (y)$, $x,y \in X$. Then,
for every function $f : X \to \rr$,
$$ {\cal E} (f,f) = \int_X f (- \L f)  d\mu = {\rm Var }_\mu (f)  $$
so that the spectral gap $\lambda =1$. The logarithmic Sobolev constant is known to be
$$ \rho = { 2 (p-q) \over  \log p - \log q} \quad (= 1 \; \; {\hbox {if}} \; \; p=q). \eqno (12)$$
The product chain on the discrete cube $X = \{-1,+1\}^N$ with the product probability measure
$\mu  = (p \delta_{+1} + q \delta_{-1})^{\otimes N}$ and generator
$\L = \sum_{i=1}^n \L_i$ is associated to the  Dirichlet form
$$ {\cal E} (f,f) = \int _X \sum_{i=1}^N f (-\L_i f)  d\mu  = pq \int _X \sum_{i=1}^N  | D_i f |^2 d\mu $$
where $D_i f$ is defined in (2). By the previous product property, it admits 1 as spectral gap
and $\rho $ given by (12) as logarithmic Sobolev constant. In its hypercontractive
formulation, the case $p= q$ is the content of the Bonami-Beckner inequality [Bo], [Be].

As mentioned before, M. Talagrand [T] used thus hypercontractivity on the discrete cube
$\{-1,+1\}^N$ equipped with the product
measure $\mu  = (p \delta_{+1} + q \delta_{-1})^{\otimes N}$
to prove that for any function $ f :  \{-1,+1\}^N \to \rr$,
$$ {\rm Var}_\mu  (f) \leq  { C pq (\log p - \log q) \over p - q} \, 
     \sum_{i=1}^N {  {\| D_i f\| }_2^2 
          \over 1+ \log \big ( {\| D_i f\| }_2 / 2\, \sqrt {pq} \, {\| D_i f\| }_1\big ) }   \eqno (13) $$
for some numerical constant $C >0$ (this statement will be covered in Section 2 below).          
This in turn yields a version of the influence result of [K-K-L] on the biased cube.   

In the continuous setting $X = \rr^n$, the case of a quadratic potential $V$ amounts to the
Hermite or Ornstein-Uhlenbeck operator $\L = \Delta - x \cdot \nabla $ with invariant measure 
the standard Gaussian measure $d\mu (x) = (2\pi )^{-n/2} \, \eu^{- |x|^2/2} dx$. It is known
here that $\lambda = \rho =1$ independently of the dimension. (More generally,
if $ V (x) - c \, {|x|^2 \over 2}$ is convex for some $c>0$, then $\lambda \geq  \rho \geq c$.)
Actually, $\L $ may also be viewed as the sum $\sum_{i=1} ^n \L_i$ of one-dimensional
Ornstein-Uhlenbeck operators along each coordinate, and $\mu $ as the product
measure of standard normal distributions. Within this product structure, the
analogue (6) of (13) has been known for some time, and will be recalled below.

\vskip 8mm  \goodbreak

{\bf 2. Hypercontractivity and Talagrand's inequality}

\vskip 3mm

This section presents the general hypercontractive approach to Talagrand type inequalities
including the discrete cube, the Gaussian product measure and more general non-product
models. The method of proof, directly inspired from [T], has been developed recently by
R. O'Donnell and K. Wimmer [OD-W1], [OD-W2] towards non-product extensions
on suitable graphs. Besides hypercontractivity, a key feature necessary to develop the
argument is a suitable decomposition of the Dirichlet form along ``directions"
commuting with the Markov operator or its semigroup. These directions are immediate in a product
space, but do require additional structure in more general contexts.

In the previous abstract setting of a Markov semigroup $\P$ with generator $\L $,
assume thus that the associated Dirichlet form ${\cal E}$ may be
decomposed along directions $\Gamma_i$ acting on functions on $X$ as
$$ {\cal E} (f, f) = \sum_{i=1}^N \int_X \Gamma _i (f)^2 d\mu  \eqno (14)$$
in such a way that, for each $i= 1, \ldots , N$, $\Gamma _i $ commutes to $\P $
in the sense that, for some constant $ \kappa  \in \rr$, every $t\geq 0$ and every $f$ in a suitable
family of functions,
$$ \Gamma _i (P_t f ) \leq  \eu^{\kappa t} \, P_t \big ( \Gamma _i (f) \big) .  \eqno (15) $$
These properties will be clearly illustrated on the main examples of interest below,
with in particular explicit descriptions of the classes of functions for which (14) and (15) may hold.

We first present the Talagrand inequality in this context. The proof is the
prototype of the hypercontractive argument used throughout this note and
applied to various examples.

\vskip 2mm

{\bf Theorem 1.} {\sl In the preceding setting, assume that $(\L , \mu )$ is hypercontractive
with constant $\rho >0$ and that (14) and (15) hold. Then, for any function $f$ in $\L^2(\mu)$,
$$ {\rm Var}_\mu (f) \leq  C(\rho , \kappa ) \sum_{i=1}^N 
       {{\| \Gamma  _i f\| }_2^2 \over  1 + \log ( {\| \Gamma _i f\|}_2 / {\| \Gamma  _i f\|}_1 ) }  $$
where $C(\rho , \kappa ) = 4 \, \eu ^{(1 + (\kappa /\rho ))^+} \! / \rho $.}
\vskip 2mm

{\it Proof.}    
The starting point is the variance representation along the semigroup $\P$
of a function $f $ in the $\L^2(\mu )$-domain of the semigroup as
$$  {\rm Var}_\mu (f) 
       = - \int_0^\infty \bigg ({d \over dt} \int_X (P_t f)^2 d \mu \bigg ) dt 
       =  - 2 \int_0^\infty \bigg ( \int_X P_t f \, \L P_t f d\mu  \bigg ) dt . $$
The time integral has to be handled both for the large and small values. For the large values
of $t$, we make use of the exponential decay provided by the spectral gap in the form
of (10) to get that, with $T =  1/2\rho  $ for example since $\rho  \leq \lambda $,
$$ {\rm Var}_\mu (f)  \leq 2 \, \big [ {\| f\| }^2_2 -  {\|  P_T f\| }^2_2 \big ] . $$
We are thus left with the variance representation of
$$   {\| f\| }^2_2 -  {\|  P_T f\| }^2_2  
    = - 2 \int_0^T \bigg ( \int_X P_t f \, \L P_t f d\mu  \bigg ) dt 
         = 2 \int_0^T {\cal E} (P_t f, P_t f)   dt . $$
Now by the decomposition (14),         
$$   {\| f\| }^2_2 -  {\|  P_T f\| }^2_2 
        = 2 \sum_{i=1}^N \int _0^T  \bigg ( \int _X   \big(\Gamma  _i (P_t f) \big ) ^2 d\mu \bigg ) dt  . $$ 
Under the commutation assumption (15),          
$$ \int _X   \big(\Gamma  _i (P_t f) \big ) ^2 d\mu
         \leq \eu^{2\kappa t} \int _X \big (P_t  \big (\Gamma _i (f) \big ) \big )^2 d\mu . $$
Since $\P $ is hypercontractive with constant $\rho >0$, for every
$i = 1, \ldots , N$ and $t \geq 0$,
$$ \big \|   P_t  \big ( \Gamma _i (f) \big ) \big \| _2  \leq  { \big \|   \Gamma  _i (f)  \big \|} _p $$     
where $ p = p(t) = 1 + \eu^{-2\rho t} \leq 2$. After the change of variables
$ p(t) = v $, we thus reached at this point the inequality
$$  {\rm Var}_\mu (f)   \leq  {2 \, \eu ^{(1 + (\kappa /\rho ))^+} \over \rho }
       \sum_{i=1}^N \int_1^2  { \big \|   \Gamma  _i (f)  \big \|} ^2_v \, dv. 
       \eqno (16) $$
This inequality actually basically amounts to Theorem 1. Indeed, by H\"older's inequality,
$$  { \big \|   \Gamma  _i (f)  \big \|} _v   \leq 
      { \big \|   \Gamma  _i (f)  \big \|} _1^\theta   \, { \big \|   \Gamma  _i (f)  \big \|} _2^{1 -\theta } $$
where $\theta  = \theta (v) \in [0,1]$ is defined by $ {1 \over v} = {\theta \over 1}  + {1-\theta \over 2}$.
Hence
$$  \int_1^2  { \big \|   \Gamma  _i (f)  \big \|} ^2_v \, dv 
     \leq   { \big \|   \Gamma  _i (f)  \big \|} ^2_2  \int _1^2  b^{2 \theta (v) } dv $$
where $ b = {\| \Gamma  _i (f)\| }_1  / {\| \Gamma  _i (f)\| }_ 2 \leq 1$. It remains to
evaluate the latter integral with $2\theta (v) = s$,
$$ \int_1^2 b ^{2 \theta (v) } dv  \leq \int_0^2 b^s ds \leq  {  2  \over  1 + \log (1/b)}   $$
from which the conclusion follows. \cqfd

\vskip 2mm

Inequality (16) of the preceding proof may also be used towards a version
of Theorem~1 with Orlicz norms as emphasized in [T]. As in [T], let $\varphi : \rr_+ \to \rr_+$ be
convex such that $\varphi (x) = {x^2/\log (\eu+x)}$ for $x\geq 1$, and $\varphi (0) = 0$, and denote
$$  {\| g\| }_\varphi  
          = \inf \bigg \{ c > 0  \, ; \int_X \varphi \big ( |g | / c \big ) d\mu  \leq 1 \bigg \} $$
the associated Orlicz norm of a measurable function $ g : X \to \rr$.
Then, for some numerical constant $C>0$,
$$ \int _1^2 {\|g \|}_v^2 \, dv \leq  C \, {\|g\|}_\varphi ^2  \eqno (17) $$
so that (16) yields
$$  {\rm Var}_\mu (f)   \leq  {2 C \, \eu ^{(1 + (\kappa /\rho ))^+} \over \rho }
       \sum_{i=1}^N   { \big \|   \Gamma  _i (f)  \big \|} ^2_\varphi  . \eqno (18) $$
Since as pointed out in Lemma 2.5 of [T],
$$  {\|g\|}_\varphi ^2  \leq   { C  \,  {\| g\|} _2^2 \over  1 + \log ( {\| g \|}_2 / {\| g\|}_1 ) } \, , $$
we see that (18) improves upon Theorem 1.
To briefly check (17), assume by homogeneity that $ \int_X g^2 /\log (e+g) d\mu  \leq 1 $
for some non-negative function $g$. Then, setting $g_k = g \, 1_{\{ 2^{k-1} <g\leq 2^k\} }$,
$k\geq 1$, and $g_0 = g \, 1_{\{ g\leq 1\} }$, 
$$ \sum_{k \in \nnn} {1 \over k+1} \int_X  g_k^2 d\mu  \leq C_1 \eqno (19)$$
for some numerical constant $C_1>0$. Hence, since $g_k \leq 2^k$ for every $k$,
$$ \eqalign {
     \int _1^2 {\| g\|}_v^2 \, dv 
    & = \int _1^2 \bigg ( \sum_{k \in \nnn} \int_X  g_k^v d\mu  \bigg )^{2/v} dv \cr
    & \leq 4 \int _1^2 \bigg ( \sum_{k \in \nnn} 2^{-(2-v)k} \int _X  g_k^2 d\mu  \bigg )^{2/v} dv \cr
    & \leq C_2 \sum_ {k \in \nnn} \bigg(\int _1^2 (k+1)^{2/v} 2^{-2(2-v)k/v} dv \bigg )
                                         {1\over k+1}\int \! g_k^2 d\mu   \cr } $$
where we used $(19)$ as convexity weights in the last step.
Now, it is easy to check that
$$ \int _1^2 (k+1)^{2/v} 2^{-2(2-v)k/v} dv \leq  C_3 $$
uniformly in $k$ so that $  \int _1^2 {\| g\|}_v^2 \, dv  \leq C_1C_2C_3$ concluding thus
the claim.

We next illustrate the general Theorem 1 on various examples of interest.

On a probability space $(X, {\cal A}, \mu )$, consider first
the Markov operator $\L f = \int_X f d\mu  - f$ acting on integrable
functions (in other words $Kf = \int_X f d\mu$).
This operator is symmetric with respect to $\mu $ with Dirichlet form
$$ {\cal E} (f, f) = \int _X f (- \L f) d\mu  = {\rm Var}_\mu (f).$$
In particular, it has spectral gap 1.
Let now $X = X_1 \times \cdots \times X_N$ be a product space with product
probability measure $ \mu  = \mu _1 \otimes \cdots \otimes \mu _N$.
Consider the product operator $\L = \sum_{i=1}^N \L_i$ where $\L_i$ is
acting on the $i$-th coordinate of a function $f$ as $\L_i f = \int_{X_i} f d\mu_i  - f$.
The product operator $\L $ has still spectral gap 1. Its Dirichlet form is given by
$$ {\cal E} (f, f) = \sum_{i=1}^N \int_X f ( - \L_i f )  d \mu 
     =   \sum_{i=1}^N \int_X  (  \L_i f )^2  d \mu  .$$
We are therefore in the setting of a decomposition of the type (14). Moreover, it is immediately
checked that $ \L_i \, \L = \L \, \L_i$ for every $i = 1, \ldots, N$, and thus
the commutation property
(15) also holds (with $\kappa =0$). Hence Theorem 1 applies for this model
with hypercontractive constant $ \rho = \min_{1 \leq i\leq N} \rho _i  >0$.
In particular, Theorem 1 includes Talagrand's inequality (13) for the hypercube $X = \{-1, +1\}^N$ with
the product measure $\mu  = (p \delta_{+1} + q \delta_{-1})^{\otimes N}$
with hypercontractive constant given by (12), for which it is
immediately checked that, for every $r\geq 1$ and every $i = 1, \ldots , N$,
$$ \int_X | \L_i f |^r d\mu  = (pq^r + p^rq) \int_ X |D_i f|^r d\mu .$$
 
Non-product examples may be considered similarly as has been thus
emphasized recently in [OD-W1] and [OD-W2] with similar arguments.
Let for example $G$ be a finite group, and let $S$ be a symmetric set
of generators of $G$. The Cayley graph associated to $S$ is the graph with vertices
the element of $G$ and edges the couples $(g,gs)$ where $g \in G$ and $s \in S$. The transition
kernel associated to this graph is
$$ K(x,y) = {1 \over |S|} \, {\bf 1}_S (y x^{-1}), \quad x, y \in G,  $$
where $|S|$ is the cardinal of $S$. The uniform probability measure $\mu $
on $G$ is an invariant and reversible measure for $K$.
This framework includes the example of $G = {\cal S}_n$ the
symmetric group on $n$ elements
with the set of transpositions as generating set and the uniform
measure as invariant and symmetric measure.

Given such a finite Cayley graph $G$ with generator set
$S $, kernel $K$ and uniform measure $\mu $ as invariant measure, the associated
Dirichlet form may be expressed on functions $ f : G \to \rr$ in the form (14)
$$ {\cal E} (f,f) = {1 \over 2 |S|} \sum_{s \in S} \sum_{x \in G} \big [ f(sx) - f (x) \big]^2 \mu (x)
   = {1 \over 2 |S|} \sum_{s \in S} {\| D _s f \| }_2^2 $$
where for $s \in S$, $ D _s f (x) = f(sx) - f(x)$, $ x \in G$.
In order that the operators $D_s$ commute to $K$ in the sense of (15)
(with again $\kappa =0$), it is necessary to assume
that $S$ is stable by conjugacy in the sense that
$$ {\hbox {for all}} \, \, u \in S, \quad u \, S \,  u^{-1} = S $$ 
as it is the case for the set of transpositions on the symmetric group ${\cal S}^n$.
The following statement from [OD-W1] is thus an immediate
consequence of the general Theorem 1.

\vskip 2mm

{\bf Corollary 2.} {\sl Under the preceding notation and assumptions, denote by
$\rho $ the logarithmic Sobolev constant of the chain $ (K, \mu )$. Then
for every function $f$ on $G$,}
$$ {\rm Var}_\mu (f ) \leq  {2 \eu \over \rho |S| } \, 
         \sum_{s \in S} {\| D _s f\| _2^2 \over 
                         1 + \log \big ( \| D _s f\| _2 / \| D _s f\| _1 \big ) } \, . $$

\vskip 2mm

One may wonder for the significance of this Talagrand type inequality for influences.
For $A \subset G$
and $ s \in S$, define the influence $I_s(A)$ of the direction $s$  on the set $A$ by
$$ I_s(A) = \mu \big ( \{ x \in G ; x \in A, sx \notin A \} \big ).$$
As on the discrete cube, given $A \subset G$ with $\mu (A) =a $, Corollary 2 yields
the existence of $s \in S$ such that
$$ I_s(A)  \geq {1 \over C} \,   a(1-a) \rho  
         \, \log \Big ( 1 +  {1\over C\rho \, a(1-a) } \Big ) 
         \geq  {1 \over C}  \, a(1-a) \,   \rho \log \Big ( 1 +  {1\over C\rho  } \Big ) 
              \eqno (20) $$
(where $C \geq 1$ is numerical).
However, with respect to the spectral gap inequality of the chain $(K, \mu )$
$$ \lambda \, {\rm Var}_\mu (f) \leq 
         {1 \over 2 |S|} \sum_{s \in S} {\| D _s f \| }_2^2 \, , $$
we see that (20) is only of interest 
provided that $\rho \log (1 + (1/\rho )) > \! \! > \lambda $. This is the case
on the symmetric discrete cube $\{- 1, +1 \}^N$ for which, in the Cayley graph
normalization of Dirichlet forms, $\lambda = \rho  = 1/N$. On the symmetric group,
it is known that the spectral gap $\lambda $  is ${2 \over  n-1}$ whereas its logarithmic
Sobolev constant $\rho $ is of the order of $1/n \log n$ ([D-SC], [L-Y]) so that
$\rho \log (1 + (1/\rho )) $ and $\lambda $ are actually of the same order for large
$n$, and hence yield the existence of a transposition $\tau $
with influence at least only of the order of $1/n$. It is pointed out in [OD-W2]
that this result is however optimal.
The paper [OD-W1] presents examples in the more general context
of Schreier graphs for which (20) yields influences strictly better than the ones from
the spectral gap inequality.

Theorem 1 may also be illustrated on continuous models such as Gaussian measures.
While the next corollary is stated in some generality, it is already of
interest for products of one-dimensional factors and covers in particular
the example (6) of the standard Gaussian product measure.

\vskip 2mm

{\bf Corollary 3.} {\sl Let $d\mu_i (x) = \eu^{-V_i(x)}dx  $, $i = 1, \ldots, N$, on
$X_i = \rr^{n_i}$ be hypercontractive with constant $\rho_i >0$.
Let $\mu = \mu _1 \otimes \cdots \otimes \mu _N$ on
$ X = X_1 \times \cdots \times  X_N$. Assume in addition that $V''_i \geq - \kappa $,
$\kappa  \in \rr $, $i = 1, \ldots , N$.  Then, for any smooth function $f$ on $X$,
$$ {\rm Var}_\mu (f) \leq  C( \rho , \kappa  ) \sum_{i=1}^N 
       {\| \nabla  _i f\| _2^2 \over  1 + \log \big ( \| \nabla _i f\|_2 / \| \nabla  _i f\|_1 \big  ) } $$
where $ \rho = \min _{1\leq i\leq N} \rho _i$, and       
where $\nabla _i f$ denotes the gradient of $f$ in the direction $X_i$, $i = 1, \ldots, N$.}  

\vskip 2mm     

Corollary 3 again follows from Theorem 1. Indeed, the
product structure immediately allows for the decomposition (14) of the Dirichlet form
$$ {\cal E} (f, f) = \int_X |\nabla f|^2 d\mu 
           = \sum_{i=1}^N \int_X |\nabla _i f|^2 d\mu $$
along smooth functions with thus $\Gamma _i(f) = |\nabla _i f|$. On the other hand, the   
basic commutation (15) between the semigroup and the gradients $\nabla _i$
is described here as a curvature condition. Namely, whenever the Hessian $V''$ of
a smooth potential $V$ on $\rr^n$ is (uniformly) bounded below by $- \kappa  $, $ \kappa  \in \rr$,
the semigroup ${(P_t)}_{t\geq 0} $ generated by the operator $\L  = \Delta - \nabla V \cdot \nabla $
commutes to the gradient is the sense that, 
for every smooth function $f$ and every $t\geq 0$,
$$ | \nabla P_t f | \leq  \eu^{\kappa t} \, P_t \big ( | \nabla f | \big ) . \eqno (21) $$
In the product setting of Corollary 3, the semigroup $\P$ is the tensor product
of the semigroups along every coordinate so that (21) ensures that
$$ | \nabla_i P_t f | \leq  \eu^{\kappa t} \, P_t \big ( | \nabla_i f | \big )  \eqno (22) $$
along the partial gradients $\nabla _i$, $i = 1, \ldots , N$ and hence (15) holds on smooth functions.
This commutation property (with $\kappa  = -1$) is for example explicit on the integral representation
$$ P_t f(x) = \int_{\rrr^n} f \big ( \eu^{-t} x + (1 - \eu^{-2t})^{1/2} y \big ) d\mu (y),
       \quad x \in \rr ^n, \, \, t \geq 0, \eqno (23) $$
of the Ornstein-Uhlenbeck semigroup
with generator $ \L = \Delta - x \cdot \nabla $
and invariant and symmetric measure the standard Gaussian distribution.  The assumption
$ V'' \geq - \kappa $ describes a curvature property of the generator $\L$ and is linked to
Ricci curvature on Riemannian manifolds. Since only $\kappa \in \rr$ is required here, it appears
as a mild property, shared by numerous potentials such as
for example double-well potentials on the line of the form $V(x) = a x^4 - bx^2$,
$a,b>0$. Recall that the assumption $V'' \geq c > 0$ (for example the quadratic potential
with the Gaussian measure as invariant measure) actually
implies that $\mu $ satisfies a logarithmic Sobolev inequality,
and thus hypercontractivity (with constant $c$).
We refer for example to [Ba], [L1], [B-G-L]... for an account on (21) and the preceding
discussion.

Corollary 3 admits generalizations in broader settings.
Weighted measures on Riemannian manifolds with a lower bound
on the Ricci curvature may be considered similarly with the same conclusions. 
In another direction, the hypercontractive approach may be developed
in presence of suitable geometric decompositions. The next statements
deal with the example of the sphere and with
geometric decompositions of the identity in Euclidean space which are
familiar in the context of Brascamp-Lieb inequalities (see [B-CE-L-M]
for further illustrations in a Markovian framework).

A non-product example in the continuous setting is the one of the standard
sphere $ \ss^{n-1} \subset \rr^n$ ($n\geq 2$) equipped with its uniform normalized measure $\mu $.
Consider, for every $i,j = 1, \ldots , n$, $D_{ij} = x_i \partial _j - x_j \partial _i$.
These will be the directions along which the Talagrand inequality may be considered since
$$ {\cal E}(f,f) = \int_{\sss^{n-1}} f (- \Delta f)  d\mu  
           = {1\over 2} \sum_{i,j=1}^n  \int_{\sss^{n-1}} (D_{ij} f)^2  d\mu .  $$
The operators $D_{ij}$ namely commute in an essential way to the spherical Laplacian
$\Delta = {1\over 2} \sum_{i,j=1}^n D_{ij}^2$ so that (15) holds with $\kappa =0$.
Finally, the logarithmic Sobolev constant  is known to be $n-1$ [Ba], [L1], [B-G-L]....
Corollary 4 thus again follows from the general Theorem 1.

\vskip 2mm

{\bf Corollary 4.} {\sl For every smooth enough function $f : \ss^{n-1} \to \rr$,}
$$ {\rm Var}_\mu (f) \leq  {4  \eu \over n}  \sum_{i,j=1}^n
       {  {\| D_{ij} f \|} _2^2 \over  1 + \log \big ( {\| D_{ij} f \|} _2 / {\| D_{ij}f \| }_1 \big ) }  \, . $$

\vskip 2mm     

Up to the numerical constant, this inequality improves upon the Poincar\'e inequality
for $\mu $ (with constant $\lambda = n-1$). 

We turn to geometric Brascamp-Lieb decompositions.
Consider thus $E_i$, $i = 1, \ldots, m$, subspaces in $\rr^n$, and $c_i > 0$, $i = 1, \ldots, m$,
such that
$$ {\rm Id}_{\rrr^n} = \sum_{i=1}^m c_i \, Q_{E_i}   \eqno (24)$$
where $Q_{E_i}$ is the projection onto $E_i$. In particular, for every $ x \in \rr^n$,
$ |x|^2 = \sum_{i=1}^m c_i | Q_{E_i} (x) |^2$ and thus, for every smooth
function $f$ on $\rr^n$,
$$ {\cal E}(f,f) = \int_{\rrr^n} |\nabla f |^2 d\mu 
      =  \sum_{i=1}^m c_i  \! \bigg ( \int_{\rrr^n} \big | Q_{E_i}(\nabla P_t f) \big |^2 d\mu \bigg ) . $$
Furthermore, $ Q_{E_i}(\nabla P_t f) = \eu^{-t} P_t (Q_{E_i} (\nabla f)) $ which may be examplified
on the representation (23) of the Ornstein-Uhlenbeck semigroup with hypercontractive
constant 1. Theorem 1 thus yields the following conclusion.

\vskip 2mm

{\bf Corollary 5.} {\sl Under the decomposition (24), for $\mu $ the standard
Gaussian measure on $\rr^n$, and for every smooth function $f$ on $\rr^n$,}
$$ {\rm Var}_\mu (f) \leq  4 \sum_{i=1}^m c_i \, 
       {\big \| Q_{E_i} (\nabla f) \big \| _2^2 \over 
             1 + \log \big ( {\| Q_{E_i} (\nabla f)\|}_2 / {\| Q_{E_i} (\nabla f) \|}_1 \big ) }  \, . $$

\vskip 2mm

\vskip 8mm  \goodbreak

{\bf 3. Hypercontractivity and geometric influences}

\vskip 3mm

In the continuous context of the preceding section, and as
discussed in the introduction, the $\L^2$-norms of gradients in Corollary 3 are not
well-suited to the (geometric) influences of [K-M-S] which require $\L^1$-norms. 
In order to reach $\L^1$-norms through the hypercontractive argument,
a further simple interpolation trick will be necessary.

To this task, we use an additional feature of the curvature condition
$V'' \geq - \kappa $, $ \kappa  \geq 0 $, namely that the action of the semigroup $\P$ with generator
$\L = \Delta - \nabla V \cdot V$ on bounded functions yields functions with
bounded gradients. More precisely (cf. [L1], [B-G-L]...),
for every smooth function $f$ with $ |f| \leq 1$, and every $ 0 < t \leq 1/2\kappa $,
$$  | \nabla  P_t f |  \leq   { 1\over \sqrt t }  \, .  \eqno (25) $$
This property may again be illustrated in case of the Ornstein-Uhlenbeck
semigroup (22) for which, by integration by parts,
$$   \nabla  P_t f (x) 
     = {\eu^{-t} \over (1 - \eu^{-2t})^{1/2} }
         \int_{\rrr^n} y \,  f \big ( \eu^{-t}x + (1 - \eu^{-2t})^{1/2} y \big ) d\mu (y) .$$

With this additional tool, the following statement then presents the expected result.
The setting is similar to the one of Corollary 3.
Dependence on $\rho $ and $\kappa $ for the constant $C'(\rho , \kappa )$ below
may be drawn from the proof. It will of course be independent of $N$.

\vskip 2mm

{\bf Theorem 6.} {\sl Let $d\mu_i (x) = \eu^{-V_i(x)}dx  $, $i = 1, \ldots, N$,
on $X_i = \rr^{n_i}$ be hypercontractive with constant $\rho_i >0$. Let
$\mu = \mu _1 \otimes \cdots \otimes \mu _N$ on
$ X = X_1 \times \cdots \times  X_N$, and set as before
$ \rho  = \min_{1\leq i\leq N} \rho _i$. Assume in addition that $V''_i \geq - \kappa $,
$\kappa  \geq 0$, $i = 1, \ldots , N$.
Then, for some constant $C'(\rho ,\kappa ) \geq 1$
and for any smooth function $f$ on $X$ such that $|f|\leq 1$,}
$$ {\rm Var}_\mu (f) \leq  C'(\rho ,\kappa ) \sum_{i=1}^N    
       {  {\| \nabla  _i f\|} _1 \big ( 1 + \| \nabla  _i f\| _1 \big )
            \over  \big [ 1 +    \log ^+ \big (1 / \| \nabla _i f\|_1 \big )  \big ]^{1/2}  }  \, . $$

\vskip 2mm  

{\it Proof.} We follow the same line of reasoning as in the proof of Theorem 1, starting
on the basis of (10) from
$$ {\| f\| }^2_2 -  {\|  P_T f\| }^2_2 
          = 2 \sum_{i=1}^N \int _0^T \bigg ( \int_X | \nabla_i P_t f|^2 d\mu \bigg ) dt
          \leq  4 \sum_{i=1}^N \int _0^T \bigg ( \int_X | \nabla_i P_{2t} f|^2 d\mu \bigg ) dt  $$
for some $T > 0$.  By (21) along each coordinate, for each $t\geq 0$, 
$$   | \nabla_i P_{2t} f | \leq  \eu^{\kappa t} \, P_t \big ( | \nabla _i P_t f | \big ) .$$
Hence, by the hypercontractivity property as in Theorem 1,
$$ { \|   \nabla _i P_{2t} f  \|} _2 \leq  \eu^{\kappa t} \, { \|  \nabla _i P_t f \|} _p $$
where $ p = p(t) = 1 + \eu^{-2\rho t} \leq  2$.
We then proceed to the interpolation trick. Namely,
by (25) and the tensor product form of the semigroup,
$ |\nabla _i P_t f | \leq t^{-1/2}$ for $0< t \leq 1/2\kappa $, so that in this range,
$$ { \|   \nabla _i P_{2t} f  \|} _2  \leq \eu^{\kappa (1 + 1/p)t} \, t^{ -(1-1/p)/2} \,  
          { \|  \nabla _i  f \|} _1^{1/p} $$
(where we used again (22)). As a consequence, provided $ T \leq 1/2\kappa $,
$$ {\| f\| }^2_2 -  {\|  P_T f\| }^2_2   \leq 
       4 \, \eu^{4\kappa T}   \sum_{i=1}^N    {\| \nabla  _i f\| }_1  \int _0^T t^{ - (1 - 1/p(t))} 
                    {\| \nabla  _i f\| }_1^{(2/p(t)) -1}   dt .    $$
                                    
We are then left with the estimate of the latter integral that only requires elementary calculus.                  
Set $ b =  {\| \nabla  _i f\| }_1$ and $\theta (t) = {2 \over p(t)} -1 \leq 1$.
Assuming $T \leq 1$, 
$$  \int _0^T t^{ - (1 - 1/p(t))} \, b^{\theta (t) } dt 
    \leq   \int _0^T  t^{ - 1/2}  \,b^{\theta (t)} dt . $$
Distinguish between two cases. When $b \geq 1$, 
$$  \int _0^T  t^{ - 1/2} \,  b^{\theta (t)} dt \leq  b \int _0^T  t^{ - 1/2} dt  \leq   2b \sqrt T. $$    
When $ b \leq 1$, use that $ \theta (t) \geq  \rho t / 2$ for every $0 \leq t\leq 1/2\rho $. Hence,
provided $ T \leq 1/2\rho $,
$$  \int _0^T t^{ -1/2}  \,  b ^{\theta (t)} dt  \leq \int _0^T  t^{ -1/2} \,   b ^{\rho t /2} dt 
        \leq  { C \over \sqrt \rho } \cdot {1 \over   \big [1 +  \log (1/b) \big ]^{1/2} }  $$
where $C \geq 1 $ is numerical.        
Summarizing, in all cases, provided $T$ is chosen smaller than
$ \min \big (1, {1 \over 2\rho } \big )$, we have
$$  \int _0^T t^{ - (1 - 1/p(t))}  b^{\theta (t)} dt  
            \leq  { 2C \over \sqrt \rho } \cdot {1 + b \over   \big [1 +  \log^+ (1/b) \big ]^{1/2} } \, . $$ 
Choosing for example
$ T = \min \big (1, {1 \over 2\rho }, {1 \over 2\kappa } \big )$ and using (10),
Theorem 6 follows with $C'(\rho , \kappa ) = C' / \rho ^{3/2} T$ for some further numerical
constant $C'$.
If $ \kappa  \leq c \rho $, then this constant is of order $\rho ^{-1/2}$. \cqfd        

\vskip 2mm

The preceding proof may actually be adapted to interpolate between Corollary 3 and
Theorem 6 as
$$ {\rm Var}_\mu (f) \leq   C \sum_{i=1}^N 
       {  {\| \nabla  _i f\|} _q^q  \, \big ( 1 +   {\| \nabla  _i f\| }_1^2 / {\| \nabla _i f\|}_q^q    \big ) 
         \over  \big [ 1 +  \log^+ \big ( {\| \nabla _i f\|}_q^q / {\| \nabla  _i f\|}_1^2\big ) \big ]^{q/2}  }  $$
for any smooth function $f$ on $X$ such that $ |f|\leq 1$, and any $1 \leq q \leq 2 $
(where $C$ depends on $\rho $, $\kappa $ and $q$).
              
As announced in the introduction, the
conclusion of Theorem 6 may be interpreted in terms of influences. Namely,
for $f = {\bf 1}_A$ (or some smooth approximation), define ${\| \nabla  _i f\| }_1$  as the
geometric influence $I_i(A)$ of the $i$-th coordinate on the set $A$. In other words, $I_i(A)$
is the surface measure of the section of $A$ along the fiber of $x \in X = X_1 \times \cdots \times X_N$
in the $i$-th direction, $1 \leq i\leq N$,
averaged over the remaining coordinates (see [K-M-S]). Then Theorem 6 yields that
$$ \mu (A) \big ( 1 - \mu (A) \big) \leq   C(\rho , \kappa ) \sum_{i=1}^N 
       {   I_i(A)  \big ( 1 + I_i(A) \big )  \over 
                \big [ 1 +   \log^+ \big (1/ I_i(A) \big )  \big ]^{1/2}  }  \, . $$
Proceeding as in the introduction for influences on the cube, the following consequence holds.

\vskip 2mm

{\bf Corollary 7.} {\sl In the setting of Theorem 6, for any Borel set $A$ in $X$
with $\mu (A) = a$, there is a coordinate $i$, $1 \leq i\leq N$, such that
$$ I_i(A) \geq  {   a (1-a)\over  C N } 
          \, \bigg ( \log {  N \over   a (1-a) } \bigg ) ^{1/2} 
                 \geq     {  a (1-a)  (\log N) ^{1/2} \over C N }   $$
where $C$ only depends on $\rho $ and $\kappa $. }              

\vskip 2mm  

It is worthwhile mentioning that when $N=1$, $I_1(A)$ corresponds to the
surface measure (Minkowski content)
$$ \mu ^+(A) = \liminf_{\varepsilon \to 0} {1\over \varepsilon } \, 
          \big [ \mu (A_\varepsilon ) - \mu (A) \big ] $$
of $A \subset \rr^{n_1}$,
so that Corollary 7 contains the quantitative form of the isoperimetric inequality for Gaussian measures
$$ \mu ^+(A)  \geq  {1 \over C}  \,  a(1-a)
          \, \bigg ( \log {  1 \over    a (1-a) } \bigg ) ^{1/2} . $$
Recall indeed (cf. e.g. [L1-2]) that the Gaussian isoperimetric
inequality indicates that $\mu ^+(A) \geq \varphi \circ \Phi ^{-1} (a)$
($ a = \mu (A)$) where
$\varphi (x) = (2\pi )^{-1/2} \, \eu^{-x^2/2}$, $x\in \rr$, $\Phi (t) = \int _{-\infty}^t \varphi (x) dx$,
$t \in \rr$, and that $\varphi \circ \Phi ^{-1}(u) \sim u (2 \log {1\over u})^{1/2}$ as $u \to 0$.
This conclusion, for hypercontractive log-concave measures,
was established previously in [B-L]. See [Mi1-2] for recent improvements in this regard.
   
Theorem 6 admits also generalizations in broader settings
such as weighted measures on Riemannian manifolds with a lower bound on the Ricci curvature
(this ensures that both (21) and (25) hold).

Besides the Gaussian measure, N. Keller, E. Mossel and A. Sen [K-M-S] also investigate with
isoperimetric tools products of one-dimensional distributions of the type
$c_\alpha \eu^{-|x|^\alpha }dx$, $1 < \alpha < \infty$, for which they produce
influences at least of the order of $ { (\log N)^{\beta/2} \over N}$ where
$ \beta  = 2 (1 - {1 \over \alpha })$ ($ \alpha = 2$ corresponding to the Gaussian case).
The proof of Theorem 6 may be adapted to cover this result but only seemingly for $1 < \alpha < 2$.
Convexity of the potentials $|x|^\alpha $ ensures (21) and (25).
When $1 < \alpha  < 2$, measures $ c_ \alpha \eu^{- |x|^\alpha } dx$ are
not hypercontractive. Nevertheless,
the hypercontractive theorems in Orlicz norms of [B-C-R] still indicate that
the semigroup $\P$ generated by the potential $|x|^\alpha $ is such that, for every
bounded function $g$ with ${\| g\| }_\infty =1$ and every $0 \leq t \leq 1$,
$$ {\| P_t g \| }_2^2 \leq  C \, {\| g\| }_1 
           \exp \Big ( -c \, t \log^\beta  \big ( 1 + (1 /\| g\| _1) \big ) \Big )  \eqno (26) $$
for $ \beta > 0$ and some constants $C,c>0$, and similarly for the product
semigroup with constants independent of $N$.
The hypercontractive step in the proof of Theorem 6 is then modified into
$$ \big \|   | \nabla _i P_{2t} f | \big \| _2 ^2
     \leq  C {\| \nabla _i f\| }_1 \int_0^1 t^{-1/2} 
              \exp \Big ( -ct \log^\beta \big (1 + (1/ {\| \nabla _i f\| }_1) \big ) \Big ) dt .$$
As a consequence, for any smooth $f$ with $|f|\leq 1$,
$$ {\rm Var}_\mu (f) \leq  C \sum_{i=1}^N   
       {  {\| \nabla  _i f\|} _1 \big (  1 + \| \nabla  _i f\| _1   \big ) 
          \over  \big [ 1 +    \log^+ \big ( 1/\| \nabla _i f\|_1 \big ) \big ]^{\beta /2} }  \, .
             \eqno (27) $$       
We thus conclude to the influence result of [K-M-S] in this range. When $\alpha >2$
($\beta \in (1,2)$), the potentials are hypercontractive in the usual sense so
that the preceding proofs
yield (27) but only for $\beta = 1$. We do not know how to reach the exponent
$\beta /2$ in this case by the hypercontractive argument.

We conclude this note by the $\L^1$ versions of Corollaries 4 and 5.
In the case of the sphere,
the proof is identical to the one of Theorem 6 provided one uses
that $  |D_{ij} f | \leq |\nabla f |$ which ensures that $  |D_{ij} P_t f | \leq 1 / \sqrt t $.
The behavior of the constant is drawn from the proof of Theorem 6.

\vskip 2mm

{\bf Theorem 8.} {\sl For every smooth enough function $f : \ss^{n-1} \to \rr$
such that $|f|\leq 1$,}
$$ {\rm Var}_\mu (f) \leq  {C \over \sqrt n}  \sum_{i,j=1}^n
       {  {\| D_{ij}  f \|} _1 \big ( 1 + {\| D_{ij} f \|} _1 \big )
            \over \big [  1 + \log ^+ \big ( 1/ {\| D_{ij} f \|} _1 \big ) \big ]^{1/2} }  \, . $$

\vskip 2mm  

Application to geometric influences $I_{ij} (A)$ as the limit of ${\| D_{ij} f \|} _1$
as $f$ approaches the characteristic function of the set $A$ may be drawn as in the previous
corresponding statements. From a geometric perspective, $I_{ij} (A)$ can be viewed
as the average over $x$ of the boundary of the section of $A$ in
the $2$-plane $x + {\rm span} (e_i, e_j)$. We do not know if the
order $n^{-1/2}$ of the constant in Theorem 8 is optimal.

As announced, the last statement is the $\L^1$-version of the 
geometric decompositions of Corollary 5 which seems again of interest for influences.
Under the corresponding commutation properties, the proof is developed similarly.

\vskip 2mm

{\bf Proposition 9.} {\sl Under the decomposition (24), for $\mu $ the standard
Gaussian measure on $\rr^n$ and for every smooth function $f$ on $\rr^n$
such that $|f|\leq 1$,
$$ {\rm Var}_\mu (f) \leq  C \sum_{i=1}^m c_i  \,   
       { {\| Q_{E_i} (\nabla f)  \|} _1 \big ( 1 + {\| Q_{E_i} (\nabla f)\|}_1 \big )   \over
              \big [ 1 +    \log^+ \big (1 / {\| Q_{E_i} (\nabla f)\|}_1 \big )\big ]^{1/2}  }   $$
where $C>0$ is numerical.}
   
\vskip 2mm     

Let us illustrate the last statement on a simple decomposition. As in the 
Loomis-Whitney inequality, consider the decomposition
$$ {\rm Id}_{\rrr^n} = \sum_{i=1}^n {1 \over n-1} \,  Q_{E_i}   $$
with $E_i = {e_i}^\perp$, $i = 1, \ldots , n$, $(e_1, \ldots , e_n)$ orthonormal basis.
Proposition 9 applied to $ f = {\bf 1}_A$ for a Borel set $ A$ in $\rr^n$ with
$\mu (A) = a$ then shows that there is a coordinate $i$, $ 1 \leq i \leq n$, such that
$$   {\big \| Q_{E_i} (\nabla f) \big \|} _1 
        \geq {1 \over C} \, a (1-a) \bigg ( \log { 1 \over  a(1-a)} \bigg)^{1/2} $$
for some constant $C>0$.
Now, $ {\| Q_{E_i} (\nabla f)\|} _1$ may be interpreted as the boundary measure of the
hyperplane section
$$A^{ x \cdot e_i } 
   = \big \{ ( x \cdot e_1, \ldots, x \cdot e_{i-1}, x \cdot e_{i+1}, \ldots, x \cdot e_n) ;
         \, ( x \cdot e_1, \ldots, x \cdot e_i, \ldots, x \cdot e_n) \in A \big \}  $$
along the coordinate $x \cdot e_i  \in \rr$ averaged over
the standard Gaussian measure.
By Fubini's theorem, there is $x \cdot e_i  \in \rr$ (or even a set
with measure as close to 1 as possible) such that
$$ \mu ^+ (A^{x \cdot e_i }) 
          \geq {1 \over C} \,  \, a (1-a) \bigg ( \log { 1 \over  a(1-a)} \bigg)^{1/2} . \eqno (28) $$
The interesting point here is that $a$ is the full measure of $A$. Indeed, recall          
that the isoperimetric inequality for $\mu $ indicates that  
$ \mu ^+ (A) \geq  \varphi \circ \Phi ^{-1} (a)$, hence a quantitative lower bound
for $ \mu ^+ (A)$
of the same form as (28). When $A$ is a half-space in $\rr^n$, thus extremal set for the
isoperimetric problem and satisfying $ \mu ^+ (A) = \varphi \circ \Phi ^{-1} (a)$, it is easy to see 
that there is indeed a coordinate $x \cdot e_i $ such that $ A^{x \cdot e_i } $ is
again a half-space in the lower-dimensional space. The preceding (28)
therefore extends this property to all sets.

\vskip 8mm

{\it Acknowledgement. We thank F. Barthe and P. Cattiaux for their help with the bound (26)
and R. Rossignol for pointing out to us the references [OD-W1] and [OD-W2].}

\vskip 10mm

\goodbreak

\baselineskip=10pt

\font\pc=cmcsc10 \rm 

\centerline {\pc References}

\vskip 0,3cm

\font\pc=cmcsc8 \rm 

\font\eightrm =cmr8 

{\eightrm 

\ref Aal|{\pc C. An\'e, S. Blach\`ere, D. Chafa\"\i, P. Foug\`eres, I. Gentil, F. Malrieu, C. Roberto,
G. Scheffer.} Sur les in\'egalit\'es de Sobolev
logarithmiques. Panoramas et Synth\`eses, vol. 10. Soc. Math. de France (2000)|

\ref Ba|{\pc D. Bakry.} L'hypercontractivit\'e et son utilisation en th\'eorie des semigroupes.
Ecole d'Et\'e de Probabilit\'es de Saint-Flour. Lecture Notes in Math. 1581, 1--114 (1994).
Springer|

\ref B-L|{\pc D. Bakry, M. Ledoux.} L\'evy-Gromov's isoperimetric inequality
for an infinite dimensional diffusion generator. Invent. math. 123, 259--281 (1996)|

\ref B-G-L|{\pc D. Bakry, I. Gentil, M. Ledoux.} Forthcoming monograph (2011)|

\ref B-C-R|{\pc F. Barthe, P. Cattiaux, C. Roberto.} Interpolated inequalities between exponential and gaussian
Orlicz hypercontractivity and isoperimetry. Revista Mat. Iberoamericana 22, 993--1067 (2006)|

\ref B-CE-L-M|{\pc F. Barthe, D. Cordero-Erausquin, M. Ledoux, B. Maurey.}
Correlation and Brascamp-Lieb inequalities for Markov semigroups (2009).
To appear in Int. Math. Res. Notices|

\ref Be|{\pc W. Beckner.} Inequalities in Fourier analysis. Ann. of Math. 102, 159Ð182 (1975)|

\ref B-H|{\pc S. Bobkov, C. Houdr\'e.} A converse Gaussian Poincar\'e-type inequality for convex functions. Statist. Probab. Lett. 44, 281--290 (1999)|

\ref Bo|{\pc A. Bonami.} \'Etude des coefficients de Fourier des fonctions de Lp(G). Ann. Inst. Fourier
20, 335Ð402 (1971)|

\ref D-SC|{\pc P. Diaconis, L. Saloff-Coste}. Logarithmic Sobolev inequalities for finite
Markov chains. Ann. Appl. Prob. 6, 695--750 (1996)|

\ref F-S|{\pc D. Falik, A. Samorodnitsky.} Edge-isoperimetric inequalities and influences.
Comb. Probab. Comp. 16, 693--712 (2007)|

\ref K-K-L|{\pc J. Kahn, G. Kalai, N. Linial.} The influence of variables on boolean functions.
29th Symposium on the Foundations of Computer Science, White Planes, 68-80 (1988)|

\ref K-M-S|{\pc N. Keller, E. Mossel, A. Sen.} Geometric influences (2010)|

\ref L1|{\pc M. Ledoux}. The geometry of Markov diffusion generators.
Ann. Fac. Sci. Toulouse IX, 305--366 (2000)|

\ref L2|{\pc M. Ledoux}. The concentration of measure phenomenon.
Math. Surveys and Monographs 89. Amer. Math. Soc. (2001)|

\ref L-Y|{\pc T. Y. Lee, H.-T. Yau.} Logarithmic Sobolev inequality for some models of random walks.
Ann. Probab. 26, 1855--1873 (1998)| 

\ref M1|{\pc E. Milman}. On the role of convexity in isoperimetry, spectral gap and concentration.
Invent. Math. 177, 1--43 (2009)|

\ref M2|{\pc E. Milman}. Isoperimetric and concentration inequalities - Equivalence
under curvature lower bound. Duke Math. J. 154,  207--239 (2010)|

\ref OD-W1|{\pc R. O'Donnell, K. Wimmer}. KKL, Kruskal-Katona, and monotone nets.
50th Annual IEEE Symposium on Foundations of Computer Science (FOCS 2009), 725--734,
IEEE Computer Soc., Los Alamitos, CA (2009)|

\ref OD-W2|{\pc R. O'Donnell, K. Wimmer}. Sharpness of KKL on Schreier graphs (2011)|

\ref Ros|{\pc R. Rossignol.} Threshold for monotone symmetric properties through a
logarithmic Sobolev inequality. Ann. Probab. 34, 1707--1725 (2006)|

\ref Roy|{\pc G. Royer.} An initiation to logarithmic Sobolev inequalities.
Translated from the 1999 French original SMF/AMS Texts and Monographs, 14.
Amer. Math. Soc. and Soc. Math. de France (2007)|

\ref T|{\pc  M. Talagrand.} On Russo's approximate zero-one law. Ann. Probab. 22,
1576--1587 (1994)|

}

\vskip 10mm

\parskip 1,2mm

\font\eightrm =cmr8  {\eightrm 

{\pc D. C.-E.:  Institut de Math\'ematiques de Jussieu,
Universit\'e Pierre et Marie Curie (Paris 6), 4, place Jussieu
75252 Paris Cedex 05, France,   
\font\eighttt =cmtt8  {\eighttt  cordero@math.jussieu.fr}}

{\pc M. L.: Institut de Math\'ematiques de Toulouse, Universit\'e de Toulouse,
31062 Toulouse, France, and Institut Universitaire de France,   
\font\eighttt =cmtt8  {\eighttt  ledoux@math.univ-toulouse.fr}} 

}

\bye